\documentclass{Poitiers}
\usepackage{amsmath}
\usepackage{epsfig} 
\usepackage{graphics} 

\def\currentvolume{X}

    \def\ppages{1--10}

\theoremstyle{definition}

\title[Theoretical optimization of finite difference schemes]
      {Theoretical optimization of finite difference schemes}

\author[Claire David, Pierre Sagaut ]{}

\email{david@lmm.jussieu.fr}

\begin{document}
\maketitle

\centerline{\scshape Claire David \footnotemark[1], Pierre Sagaut
\footnotemark[1] }
\medskip
{\footnotesize \centerline{\footnotemark[1] Universit\'e Pierre et
Marie Curie-Paris 6}
  \centerline{Laboratoire de Mod\'elisation en M\'ecanique, UMR CNRS 7607}
   \centerline{Bo\^ite courrier $n^0 162$, 4 place Jussieu, 75252 Paris, cedex 05,
France}

} 

\medskip

\begin{abstract}
  The aim of this work is to develop general optimization methods for finite difference schemes
  used to approximate linear differential equations. The specific
  case of the transport equation is exposed. In particular, the minimization of the numerical error is taken
  into account. The theoretical study of a related linear algebraic
  problem gives general results which can lead to the determination of the optimal scheme.
\end{abstract}

\section{Introduction: Scheme classes}
\label{sec:intro}

\indent Finite difference schemes used to approximate linear
differential equations induce numerical errors, that are generally
difficult to predict. The usual process consists in testing various
schemes for
more and more refined time and space steps.\\
We here propose a completely different approach, which consists in
determining the minimum norm error of a given finite difference
scheme. This process has the advantage of avoiding scheme
convergence tests. Moreover, it can explain error jumps that
often occur in such approximations.\\
\\

\noindent Consider the transport equation:
\begin{equation}
\label{transp} \frac{\partial u}{\partial t}+c\,\frac{\partial
u}{\partial x}=0 \end{equation}

\noindent with the initial condition $u(x,t=0)=u_0(x)$.

\bigskip

\noindent A finite difference scheme for this equation can be
written under the form:
\begin{equation} \label{scheme} {{{{{\alpha \, u}}_i}}^{n+1}}+
   {{{{{\beta \,u}}_i}}^{n}}
    +{{{{{\gamma \,u}}_i}}^{n-1}}
      +\delta \,{{{u_{i+1}}}^n}+{{{{{\varepsilon \, u}}_{i-1}}}^n}
            +{{{{{\zeta \,u}}_{i+1}}}^{n+1}}
              +{{{{{\eta \,u}}_{i-1}}}^{n-1}}+{{{{{\theta \,u}}_{i-1}}}^{n+1}}+\vartheta \,{{ u}_{i+1}}^{n-1} =0
              \end{equation}

\noindent where:
\begin{equation}
{u_l}^m=u\,(l\,h, m\,\tau)
\end{equation}
\noindent  $l\, \in \, \{i-1,\, i, \, i+1\}$, $m \, \in \, \{n-1,\,
n, \, n+1\}$, $j=0, \, ..., \, n_x$, $n=0, \, ..., \, n_t$, $h$,
$\tau$ denoting respectively the mesh size and time step.\\
The Courant-Friedrichs-Lewy number ($cfl$) is defined as $\sigma = c \,\tau / h$ .\\
\\

A numerical scheme is  specified by selecting appropriate values of
the coefficients  $\alpha$, $\beta$, $\gamma$, $\delta$,
$\varepsilon$, $\zeta$, $\eta$, $\theta$ and  $\vartheta$ in
equation (\ref{scheme}). Values corresponding to numerical schemes
retained for the present works are given in Table \ref{SchemeTable}.

\bigskip

\begin{table}
\caption{Numerical scheme coefficient.}
\begin{center} \footnotesize
{\begin{tabular}{cccccccccc} \hline
Name & $  \alpha $ & $ \beta$ & $\gamma$ & $\delta $  & $\epsilon$ & $\zeta$ & $\eta$ & $\theta$ & $\vartheta$ \\
\hline
Leapfrog & $\frac{1}{2 \tau} $ & 0 &  $\frac{-1}{2 \tau} $ & $  \frac{c}{ 2 h}  $ & $   \frac{-c}{ 2 h}  $ & 0 & 0 & 0  & 0 \\
Lax & $\frac{1}{ \tau} $ &  0 & 0 & $  \frac{-1}{ 2 \tau} + \frac{c}{ 2 h } $ & $  \frac{-1}{2 \tau}  - \frac{c}{ 2 h} $ & 0 & 0 & 0  & 0 \\
Lax-Wendroff & $\frac{1}{ \tau} $ & $  \frac{-1}{ \tau}  + \frac{c^2 \tau}{ h ^2} $ & 0 & $ \frac{( 1- \sigma ) c}{ 2 h} $&  $  \frac{-( 1+ \sigma ) c}{ 2 h } $ & 0 & 0 & 0 & 0 \\
Crank-Nicolson & $ \frac{1}{ \tau}  + \frac{c}{ h ^2} $ &  $
\frac{-1}{ \tau}  + \frac{c}{ h ^2} $ & 0 & $ \frac{-c}{ h ^2} $ & $
\frac{-c}{ h ^2 } $ & 0 & $ \frac{-c}{ h ^2} $ & $ \frac{-c}{ h ^2}
$ & 0
\end{tabular}}
\end{center}
\label{SchemeTable}
\end{table}

The number of time steps will be denoted $n_t$, the number of space
steps, $n_x$. In general, $n_t\gg n_x$.\\
\\

\bigskip

The paper is organized as follows. The equivalent matrix equation is
exposed in section \ref{Sylv}. Scheme optimization is presented in
section \ref{Opt}.

\section{The Sylvester equation}
\label{Sylv}

\subsection{Matricial form of the finite differences problem}

\noindent Let us introduce the rectangular matrix defined by:
\begin{equation}
U=[{{{u_i}}^n}{]_{\, 1\leq i\leq {n_x-1}, \, 1\leq n\leq {n_t} \, }}
\end{equation}
\\
The problem (\ref{scheme}) can be written under the following
matricial form:
\begin{equation}{M_1}\,U +U\,M_2+{\mathcal{L}}(U)=M_0
\end{equation}

\noindent where $M_1$, $M_2$ and $M_0$ are square matrices
respectively $n_x-1$ by $n_x-1$, $n_t$ by $n_t$, given by:
\begin{equation}
\begin{array}{ccc}
{M_1}= \left (
\begin{array}{ccccc}
\beta & \delta & 0 & \ldots &  0 \\
\varepsilon& \beta & \ddots & \ddots &  \vdots \\
0 &  \ddots & \ddots & \ddots &  0\\
\vdots &  \ddots & \ddots & \beta &  \delta\\
0 &  \ldots & 0 & \varepsilon &  \beta\\
\end{array} \right ) &
 & {M_2}= \left (
\begin{array}{ccccc}
0 & \gamma& 0 & \ldots &  0 \\
\alpha & 0 & \ddots & \ddots &  \vdots \\
0 &  \ddots & \ddots & \ddots &  0\\
\vdots &  \ddots & \ddots & \ddots &  \gamma\\
0 &  \ldots & 0 & \alpha &  0\\
\end{array} \right )
\end{array}
  \end{equation}

\bigskip

\bigskip
  \begin{equation}
  \label{M0}
\scriptsize{{M_0}= \left (
\begin{array}{ccccc}
    -\gamma \,u_1^0
      -\varepsilon \,u_{0}^1-\eta\,u_0^0-\theta\,u_{0}^{2}-\vartheta\,u_{2}^{0}
            &   -\varepsilon \,u_{0}^2 -\eta\,u_0^1-\theta\,u_{0}^{3}& \ldots & \ldots   &
            -\varepsilon \,u_{0}^{n_t}-\eta\,u_{0}^{n_t-1} \\
 -\gamma \,u_2^0-\eta\,u_{1}^{0}-\vartheta\,u_{3}^{0}
       &  0
 & \ldots  & \ldots    &  0 \\
\vdots  &  \vdots & \vdots & \vdots &  \vdots\\
 -\gamma \,u_{n_x-2}^0-\eta\,u_{n_x-2}^{0}-\vartheta\,u_{n_x-1}^{0}
       &  0
 & \ldots  & \ldots   &  0 \\
    -\gamma \,u_{n_x-1}^0-\delta \,u_{n_x}^{1}-\eta\,u_{n_x-2}^{0}-\zeta\,u_{n_x}^{2}-\vartheta\,u_{n_x}^{0}
      &   -\delta \,u_{n_x}^{2} -\zeta\,u_{n_x}^{3}-\vartheta\,u_{n_x}^{1}& \ldots
       & \ldots &  -\delta \,u_{n_x}^{n_t}-\vartheta\,u_{n_x}^{n_t-1}\\
\end{array} \right )}
  \end{equation}

\bigskip

\bigskip
\noindent and where ${\mathcal {L}}$  is a linear matricial operator
which can be written as: \begin{equation} {\mathcal {L}}={\mathcal
{L}}_1+{\mathcal {L}}_2+{\mathcal {L}}_3+{\mathcal {L}}_4
\end{equation}
\noindent where ${\mathcal {L}}_1$, ${\mathcal {L}}_2$, ${\mathcal
{L}}_3$ and ${\mathcal {L}}_4$ are given by:
\begin{equation}
\begin{array}{ccc}
 {\mathcal {L}}_1(U)  = \zeta  \left (
\begin{array}{ccccc}
u_2^2 & u_2^3 & \ldots & u_2^{n_t} &  0\\
u_3^2 & u_3^3 & \ldots & \vdots &   \vdots \\
\vdots &  \vdots & \ddots & \vdots &  \vdots\\
u_{n_x-1}^2 & u_{n_x-1}^3  & \ldots &  u_{n_x-1}^{n_t}& 0  \\
0 & 0  & \ldots &  0& 0  \\
\end{array}
 \right )
  &
  &  {\mathcal {L}}_2(U) = \eta \left (
\begin{array}{ccccc}
0 & 0  & \ldots &  0& 0 \\
0 & u_1^1  & u_1^2 & \ldots &  u_1^{n_t-1} \\
0 & u_1^0  & u_1^1 & \ldots &  u_2^{n_t-1} \\
\vdots &  \vdots & \vdots &  \ddots & \vdots \\
0 & u_{n_x-2}^1  & u_{n_x-2}^2 & \ldots &  u_{n_x-2}^{n_t-1} \\
\end{array}
 \right )
\end{array}
  \end{equation}

\begin{equation}
\begin{array}{ccc}
{\mathcal {L}}_3(U) = \theta \left (
\begin{array}{ccccc}
0 & \ldots &  \ldots &    \ldots &  0 \\
u_1^2 & u_1^3  & \ldots &  u_1^{n_t}& 0 \\
u_2^2 & u_2^3  & \ldots &  u_2^{n_t}& 0 \\
\vdots &  \vdots & \vdots &  \vdots & \vdots \\
u_{n_x-2}^2 & u_{n_x-2}^3  & \ldots &  u_{n_x-2}^{n_t}& 0 \\
\end{array}
 \right )
  &
  &  {\mathcal {L}}_4(U) = \vartheta \left (
\begin{array}{ccccc}
0 & u_2^1 & u_2^2 & \ldots  &  u_2^{n_t-1} \\
0 & u_3^1 & u_3^2 & \ldots  &  u_3^{n_t-1} \\
\vdots &  \vdots & \ddots  & \ddots &  \vdots\\
0 & u_{n_x-1}^1& \ldots & \ldots  &  u_{n_x-1}^{n_t-1} \\
0 &  0 & \ldots &   \ldots &  0\\
\end{array}
 \right )
 \end{array}
 \end{equation}
\bigskip

\bigskip

\noindent The second member matrix $M_0$ bears the initial
conditions, given for the specific value $n=0$, which correspond to
the initialization process when computing loops, and the boundary
conditions, given for the specific values $i=0$, $i=n_x$.
\bigskip

\noindent Denote by $u_{exact}$ the exact solution of (\ref{transp}).\\
\noindent The $U_{exact}$ corresponding matrix will be:

\begin{equation}
U_{exact}=[{{{U_{{exact}_i}}}^n}{]_{\, 0\leq i\leq {n_x-1},\, 0\leq
n\leq {n_t}\, }} \end{equation} where:

\begin{equation}
{U_{exact}}_i^n=U_{exact}(x_i,t_n)
 \end{equation}

\noindent with $x_i=i \; h$, $t_n=n \; \tau$.\\
\\

\noindent $U$ is then solution of:

   \begin{equation}
   \label{SylvCompl}
{M_1}\,U+U\,M_2+{\mathcal{L}}(U)=M_0
   \end{equation}

\noindent We will call \textit{error matrix} the matrix defined by:
 \begin{equation}
 \label{err}
E=U-U_{exact}
   \end{equation}

\noindent Let us consider the matrix $F$ defined by:
\begin{equation} F={M_1}\,U_{exact}+U_{exact}\,M_2 + {\mathcal{L}}(U_{exact})-M_0\end{equation}

\noindent The \textit{error matrix} $E$ satisfies then:

   \begin{equation}
   \label{eqmtr}
{M_1}\,E+E\,M_2+{\mathcal{L}}(E)=F
   \end{equation}

\subsection{The matrix equation}

\subsubsection{Theoretical formulation}

\noindent Minimizing the error due to the approximation induced by
the numerical scheme is equivalent to minimizing the norm of the matrices $E$ satisfying (\ref{eqmtr}).\\
\\
\noindent Since the linear matricial operator ${\mathcal{L}}$
appears only in the Crank-Nicholson scheme, we will restrain
our study to the case ${\mathcal{L}}=0$. The generalization to the case ${\mathcal{L}} \neq 0$ can be easily deduced.\\

   \noindent The problem is then the determination of the minimum norm solution
of:

   \begin{equation}
   \label{SylvErr}
{M_1}\,E+E\,M_2=F
   \end{equation}

\noindent which is a specific form of the Sylvester equation:

   \begin{equation}
   \label{SylvGen}
AX+XB=C
   \end{equation}
where $A$ and $B$ are respectively $m$ by $m$ and $n$ by $n$
matrices, $C$ and $X$, $m$ by $n$ matrices.

 \noindent The solving of the Sylvester
equation is generally based on Schur decomposition: for a given
square $n$ by $n$ matrix $A$, $n$ being an even number of the form
$n=2\,p$, there exists a unitary matrix $U$ and a upper triangular
block matrix $T$ such that:
\begin{equation}
A = U^*TU
   \end{equation}
\noindent where $U^*$ denotes the (complex) conjugate matrix of the
transposed matrix $^T U$. The diagonal blocks of the matrix $T$
correspond to the complex eigenvalues $\lambda_i$ of $A$:
\begin{equation}
T = \left (
\begin{array}{ccccc}
T_1 & 0 & \ldots  & \ldots   & 0 \\
0 & \ddots & \ddots   & \ddots &  \vdots \\
\vdots &  \ddots & T_i  & \ddots &  \vdots \\
\vdots &  \ddots & \ddots &   \ddots  & 0 \\
0 &  0 & \ldots  & 0 & T_p\\
\end{array}
 \right )
   \end{equation}
\noindent where the block matrices $T_i$, $i=1,\ ...,\, p$ are given
by:
\begin{equation}
\left (\begin{array}{cc}
\mathcal{R}e\, [\lambda_i] & \mathcal{I}m\, [\lambda_i]  \\
-\,\mathcal{I}m\, [\lambda_i] & \mathcal{R}e\, [\lambda_i]  \\
\end{array}
 \right )
   \end{equation}

\noindent $\mathcal{R}e$ being the real part of a complex number,
and $\mathcal{I}m$ the imaginary one.\\
\bigskip \noindent Due to this decomposition, the Sylvester equation
require, to be solved, that the dimensions of the matrices be even
numbers. We will therefore, in the following, restrain our study to
$n_x$ and $n_t$ being even numbers. So far, it is interesting to
note that the Schur decomposition being more stable for higher order
matrices, it perfectly fits finite differences problems.

\bigskip

\noindent Complete parametric solutions of the generalized Sylvester
equation (\ref{Sylvmod}) is given in \cite{Berman}, \cite{Gail}.\\

\noindent As for the determination of the solution Sylvester
equation, it is a major topic in control theory, and has been the
subject of numerous works (see \cite{Van}, \cite{Hearon},
\cite{Tsui}, \cite{Zhou},
\cite{Duan}, \cite{Duan2}, \cite{Kir}).\\
\noindent In \cite{Van}, the method is based on the reduction of the
he observable pair $(A,C)$ to an observer-Hessenberg  pair $(H, D)$,
$H$ being a block upper Hessenberg matrix. The reduction to the
observer-Hessenberg form $(H, D)$ is achieved by means of the
staircase algorithm
(see \cite{Bol}, ...). \\
\noindent In \cite{Zhou}, in the specific case of $B$ being a
companion form matrix, the authors propose a very neat general
complete parametric solution, which is expressed in terms of the
controllability of the matrix pair $(A,B)$, a symmetric matrix
operator, and a parametric matrix in the Hankel form.\\
We recall that a companion form, or Frobenius matrix is one of the
following kind:
 \begin{equation}
B = \left (
\begin{array}{cccccc}
0 & \ldots & \ldots  & \ldots  & 0 &  -b_{0} \\
1 & 0 & \ldots  & \ldots & 0 &  -b_{1} \\
0 &  1 & 0  & \ldots & \vdots &  \vdots \\
\vdots &  0 & \ddots & \vdots &  \vdots & \vdots \\
0 &  0 & \ldots & 1 & 0 &  -b_{p-1}\\
\end{array}
 \right )
   \end{equation}
These results can be generalized through matrix block decomposition
to a block companion form matrix, which happens to be the case of
our matrix $M_2$ in the specific case of $n_x$ and $n_t$ being even
numbers:
\begin{equation}
M_2 =  \left (
\begin{array}{ccccc}
{M_2^{B}}^1 & 0 & \ldots  & \ldots   &  0 \\
0 &  {M_2^{B}}^2 & 0 &  \ldots &  0 \\
0 &  0 & \ddots  & \ddots  &  \vdots \\
\vdots &  0 & \ddots & \ddots  & 0 \\
0 &  0 &  \ldots & 0 &  {M_2^{B}}^k\\
\end{array}
 \right )
   \end{equation}

 \noindent the ${M_2^{B}}^p$, $1\leq p \leq k$ being companion
 form matrices.\\
\\
\\
Another method is presented in \cite{Var}, where the determination
of the minimum-norm solution
of a Sylvester equation is specifically developed.\\
\noindent The accuracy and computational stability of the solutions
is examined in \cite{Deif}.
\bigskip

\subsubsection{Existence condition of the solution }

\noindent Equation (\ref{SylvGen}) has a unique
solution if and only if $A$ and $B$ have no common eigenvalues.\\
In our case, since $M_2$ is a upper triangular matrix whose diagonal
coefficients are all equal to $\alpha$, its eigenvalues are also all
equal to $\alpha$. As for the matrix $M_1$, one can easily check
that $\alpha$ does not belong to its spectra. Hence, (\ref{eqmtr})
has a unique solution, which accounts for the consistency of the
given problem.

\newpage

\section{Scheme optimization}
\label{Opt}

\noindent Advect a sinusoidal
signal\\
 \begin{equation} \label{signal} u=\text{Cos}\,[\,\frac {2\, \pi}{\lambda} \,(x-c\,t)\,]
\end{equation}
\noindent through the Lax scheme, where:
 \begin{equation} \lambda={{n_\lambda} \,dx}
\end{equation}
\noindent $n_\lambda$ denotes the number of cells per wavelength.

\bigskip

\noindent Let $n_\lambda$ remain unknown.\\

\noindent Equation (\ref{SylvErr}) can thus be normalized as:

   \begin{equation}
   \label{Sylvmod_cfl}
\overline{{M_1}}\,E+E\,\overline{M_2}=\overline{F}
   \end{equation}

\noindent where

   \begin{equation}
   \left\lbrace
   \begin{array}{rcl}
\overline{{M_1}}&=&\frac {h\,cfl}{c}\, M_1 \\
\overline{{M_2}}&=&\frac {h\,cfl}{c}\, M_2 \\
\overline{{F}}&=&\frac {h\,cfl}{c}\, F
\end{array} \right.
   \end{equation}

\noindent We deliberately choose a small value for the number of
steps: $n_t=n_x=20$, starting from the point that if the error is
minimized
 for a small value of this number, it will be the same as this
 number increases.

\noindent Figure \label{Isovaleurs} displays the $L_2$ norm of the
isovalues of the  error as a function of $n_\lambda$ (maximums are
in white, minimums in black; larger values are shown lighter).

\begin{figure}[!hbp]
    \resizebox*{0.47\columnwidth}{0.25\textheight}{\includegraphics{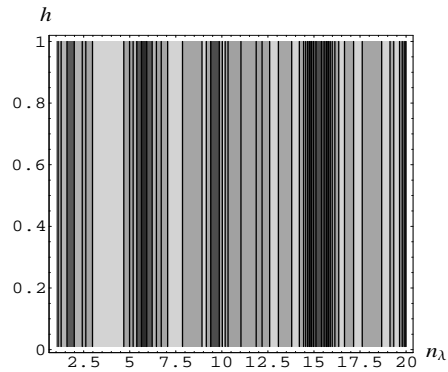}}
                \caption{Isovalues of the $L_2$ norm of the error as a
function of the number of cells per wavelength $n_\lambda$ }
                \label{Isovaleurs}
                \end{figure}

\newpage

\noindent The square value of the $L_2$ norm of the error, for two
significative values of the number of cells per wavelength $n_p$, is
displayed in Figure \ref{ErreurL2}:

\begin{figure}[!hbp]
    \resizebox*{0.87\columnwidth}{0.25\textheight}{\includegraphics{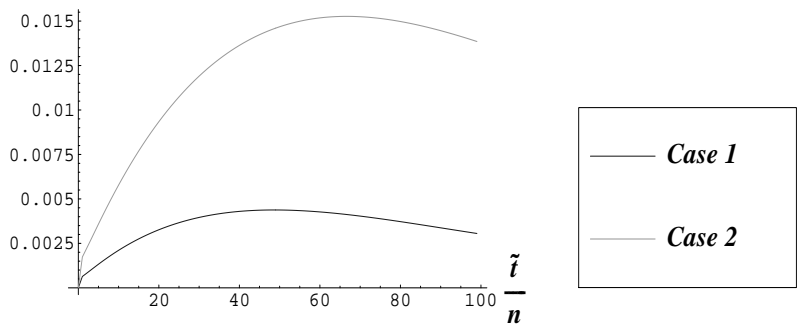}}
                \caption{square value of the $L_2$ norm  of the error as a
function of the number of cells per wavelength $n_\lambda$. Case 1:
$\lambda=9$. Case 2: $\lambda=9.8$.}
                \label{ErreurL2}
                \end{figure}

\noindent The above results ensure the faster convergence of the
error.

\section{Conclusion}

Thanks to the above results, we presently propose to optimize finite
difference problems through minimization of the symbolic expression
of the error as a  function of the
 scheme parameters.

\addcontentsline{toc}{section}{\numberline{}References}

\end{document}